\documentclass[11pt]{amsart}
\usepackage{amsmath}
\usepackage{amssymb, amsfonts, amsmath, amsthm, amscd}

\title[Projective normality of algebraic curves and its application]
      {Projective normality of algebraic curves and its application to surfaces}
\author[Seonja Kim]{Seonja Kim}
\thanks{This work was supported by the Korea Research Foundation
Grant funded by the Korean Government (R04-2003-000-10036-0 and
2005-070-C00005) for first author. This work was supported by the
Korea Research Foundation Grant funded by the Korean Government
(R14-2002-007-01000-1 and 2005-070-C00005) for second author.}
\address{Seonja Kim, \newline Department of Electronics, Chungwoon University, Chungnam, 350-701, Korea}
\email{sjkim{\char "40}chungwoon.ac.kr}
\author[Young Rock Kim]{Young Rock Kim}
\address{Young Rock Kim, \newline Department of Mathematics, Konkuk University, Seoul, 143-701, Korea}
\email{rocky777{\char "40}math.snu.ac.kr}

\def\cli{\hbox{\rm Cliff}}

\begin{document}

\newtheorem{lem}{Lemma}[section]
\newtheorem{thm}[lem]{Theorem}
\newtheorem{claim}{Claim}
\newtheorem{prop}{Theorem}
\newtheorem{Prop}[lem]{Proposition}
\newtheorem{rmk}[lem]{Remark}
\newtheorem{cor}[lem]{Corollary}
\renewcommand{\theprop}{\Alph{prop}}
\maketitle

\begin{abstract}
\noindent Let $L$ be a very ample line bundle on a smooth curve $C$
of genus $g$ with $\frac{3g+3}{2}<\deg L\le 2g-5$. Then $L$ is
normally generated if $\deg L>\max\{2g+2-4h^1(C,L),
2g-\frac{g-1}{6}-2h^1(C,L)\}$. Let $C$ be a triple covering of genus
$p$ curve $C'$ with $C\stackrel{\phi}\rightarrow C'$ and $D$ a
divisor on $C'$ with $4p<\deg D< \frac{g-1}{6}-2p$. Then
$K_C(-\phi^*D)$ becomes a very ample line bundle which is normally
generated. As an application, we characterize some smooth projective
surfaces.
\end{abstract}

\setcounter{section}{0}
\section{Introduction}

We work over the algebraically closed field of characteristic zero.
Specially the base field is the complex numbers in considering the
classification of surfaces. A smooth irreducible algebraic variety
$V$ in $\mathbb P^r$ is said to be projectively normal if the
natural morphisms $H^0(\mathbb P^r,{\mathcal O}_{\mathbb P^r}(m))\to
H^0(V,{\mathcal O}_V(m))$ are surjective for every nonnegative
integer $m$. Let $C$ be a smooth irreducible algebraic curve of
genus $g$. We say that a base point free line bundle $L$ on $C$ is
normally generated if $C$ has a projectively normal embedding via
its associated morphism $\phi_L:C\to \mathbb P(H^0(C,L))$.

Any line bundle of degree at least $2g+1$ on a smooth curve of genus
$g$ is normally generated but a line bundle of degree at most $2g$
might fail to be normally generated (\cite{KK1}, \cite{LM1},
\cite{Mum}). Green and Lazarsfeld showed a sufficient condition for
$L$ to be normally generated as follows (\cite{GL}, Theorem 1): If
$L$ is a very ample line bundle on $C$ with $\deg L\ge
2g+1-2h^1(C,L)-\cli(C)$ (and hence $h^1(C,L)\le 1$), then $L$ is
normally generated. Using this, we show that a line bundle $L$ on
$C$ with $\frac{3g+3}{2}<\deg L\le 2g-5$ is normally generated for
$\deg L>\max\{2g+2-4h^1(C,L), ~2g-\frac{g-1}{6}-2h^1(C,L)\}$. As a
corollary, if $C$ is a triple covering of genus $p$ curve $C'$ with
$C\stackrel{\phi}\rightarrow C'$ then it has a very ample
$K_C(-\phi^*D)$ which is normally generated for any divisor $D$ on
$C'$ with $4p<\deg D< \frac{g-1}{6}-2p$. It is a kind of
generalization of the result that $K_C(-rg^1_3)$ on a trigonal curve
$C$ is normally generated for $3r\le \frac{g}{2}-1$ (\cite{KK}).

As an application to nondegenerate smooth surface $S\subset \mathbb
P^r$ of degree $2\Delta-e$ with $g(H)=\Delta+f$, $\max\{\frac{e}{2},
6e-\Delta\}<f-1<\frac{\Delta-2e-6}{3}$ for some $e,f\in \mathbb
Z_{\ge 1}$, we obtain that $S$ is projectively normal with $~p_g=f$
and $-2f-e+2\le K_S^2\le \frac{(2f+e-2)^2}{2\Delta-e}$ if its
general hyperplane section $H$ is linearly normal, where
$\Delta:=\deg S-r+1$. Furthermore we characterize smooth projective
surfaces $S$ for $K_S^2=-2f-e+2$, $0$ (cf. Theorem \ref{thm5.3}).
These applications were derived by the methods in Akahori's,
Livorni's and Sommese's papers (\cite{Ak}, \cite{Liv}, \cite{So}).

We follow most notations in \cite{ACGH}, \cite{GH}, \cite{H}. Let
$C$ be a smooth irreducible projective curve of genus $g\ge 2$. The
Clifford index of $C$ is taken to be $\cli(C)=\min
\{~\cli(L)~|~h^0(C,L)\ge 2,~h^1(C,L)\ge 2~\},$ where $\cli(L)=\deg
L-2(h^0(C,L)-1)$ for a line bundle $L$ on $C$. By abuse of notation,
we sometimes use a divisor $D$ on a smooth variety $V$ instead of
${\mathcal O}_V(D)$. We also denote $H^i(V,{\mathcal O}_V(D))$ by
$H^i(V,D)$ and $h^0(V,L)-1$ by $r(L)$ for a line bundle $L$ on $V$.
We denote $K_V$ a canonical line bundle on a smooth variety $V$.

\section{Normal generation of a line bundle on a smooth curve}

Any line bundle of degree at least $2g+1$ on a smooth curve of genus
$g$ is normally generated. If the degree is at most $2g$, then there
are curves which have a non normally generated line bundle of given
degree (\cite{KK1}, \cite{LM1}, \cite{Mum}). In this section, we
investigate the normal generation of a line bundle with given degree
on a smooth curve under some condition about the speciality of the
line bundle.
\begin{thm}
Let $L$ be a very ample line bundle on a smooth curve $C$ of genus
$g$ with $\frac{3g+3}{2}<\deg L\le 2g-5$. Then $L$ is normally
generated if $\deg L>\max\{2g+2-4h^1(C,L),
2g-\frac{g-1}{6}-2h^1(C,L)\}$. \label{3.5.6}
\end{thm}
\begin{proof}

Suppose $L$ is not normally generated. Then there exists a line
bundle $A\simeq L(-R), ~R>0,$ such that \hbox{\rm (i)} $\cli(A)\le
\cli(L)$, \hbox{\rm (ii)} $\deg A\ge \frac{g-1}{2}$, \hbox{\rm
(iii)} $h^0(C,A)\ge 2$ and $h^1(C,A)\ge h^1(C,L)+2$ by the proof of
Theorem 3 in \cite{GL}. Assume $\deg K_CL^{-1}=3$, then
$|K_CL^{-1}|=g^1_3$. On the other hand, $L=K_C(-g^1_3)$ is normally
generated. 
So we may assume $\deg K_CL^{-1}\ge 4$ and then $r(K_CL^{-1})\ge 2$
since $\deg L>2g+2-4h^1(C,L)$. Set $B_1$(resp. $B_2$) is the base
locus of $K_CL^{-1}$(resp. $K_CA^{-1}$). And let
$N_1:=K_CL^{-1}(-B_1), ~N_2:=K_CA^{-1}(-B_2)$. Then $N_1\lneq N_2$
since $A\cong L(-R), ~R>0$ and $h^1(C,A)\ge h^1(C,L)+2$. Hence we
have the following diagram,

\begin{picture}(300,100)

\put(80,75){$C$}

\put(95,80){\vector(1,0){85}}

\put(90,65){\vector(2,-1){90}}

\put(186,75){$C_2$}

\put(190,70){\vector(0,-1){40}}

\put(200,45){$\pi$: projection}

\put(130,87){$\phi_{N_2}$}

\put(115,33){$\phi_{N_1}$}

\put(186,15){$C_1$}

\end{picture}

\noindent where $C_i=\phi_{N_i}(C)$.

If we set $m_i:=\deg \phi_{N_i}, ~i=1,2$, then we have $m_2|m_1$. If
$N_1$ is birationally very ample, then by Lemma 9 in \cite{KK1} and
$\deg K_CL^{-1}<\frac{g-1}{2}$ we have $r(N_1)\le \left [\frac{\deg
N_1-1}{5}\right ].$ It is a contradiction to $\deg L>2g+2-4h^1(C,L)$
that is equivalent to $\deg K_CL^{-1}<4(h^0(C,K_CL^{-1})-1)$.
Therefore $N_1$ is not birationally very ample, and then we have
$m_1\le 3$ since $\deg K_CL^{-1}<4(h^0(C,K_CL^{-1})-1)$.

Set $H_1$ be a hyperplane section of $C_1$. If $|H_1|$ on a smooth
model of $C_1$ is special, then $r(N_1)\le\frac{\deg N_1}{4}$, which
is absurd. Thus $|H_1|$ is nonspecial. If $m_1=2$, then
$$r(K_CL^{-1}(-B_1+P+Q))\ge r(K_CL^{-1}(-B_1))+1$$ for any pairs $(P,Q)$
such that $\phi_{N_1}(P)=\phi_{N_1}(Q)$ since $|H_1|$ is nonspecial.
Therefore we have $r(L(-P-Q))\ge r(L)-1$ for $(P,Q)$ such that
$\phi_{N_1}(P)=\phi_{N_1}(Q)$, which contradicts that $L$ is very
ample. Therefore we get $m_1=3$. Suppose $B_1$ is nonzero. Set $P\le
B_1$ for some $P\in C$. Consider $Q,R$ in $C$ such that
$\phi_{N_1}(P)=\phi_{N_1}(Q)=\phi_{N_1}(R)=P'$ for some $P'\in C_1$.
Since $|H_1|$ is nonspecial, we have
\begin{eqnarray*}
r(K_CL^{-1}(Q+R))&\ge& r(N_1(P+Q+R))=r(H_1+P')\\
&=&r(H_1)+1=r(K_CL^{-1})+1
\end{eqnarray*} which is a contradiction to the very ampleness of
$L$. Hence $K_CL^{-1}$ is base point free, i.e., $K_CL^{-1}=N_1$. On
the other hand, we have $m_2=1$ or 3 for $m_2| m_1$. Since
$K_CA^{-1}(-B_2)=N_2\gneq N_1=K_CL^{-1}$, we may set $N_1=N_2(-G)$
for some $G>0$.

Assume $m_2=1$, i.e. $K_CA^{-1}(-B_2)=N_2$ is birationally very
ample. On the other hand we have $r(N_2)\ge r(N_1)+\frac{\deg
G}{2}$, since $N_2(-G)\cong N_1$ and $\cli(N_2)\le
\cli(A)\le\cli(L)=\cli(N_1)$. In case $\deg N_2\ge g$ we have $n\le
\frac{2\deg N_2-g+1}{3}$ by Castelnuovo's genus bound and hence
$$\cli(L)\ge \cli(N_2)\ge \deg N_2-\frac{4\deg N_2-2g+2}{3}=\frac{2g-2-\deg N_2}{3}
\ge \frac{g-1}{6},$$ since $N_2=K_CA^{-1}(-B_2)$ and $\deg A\ge
\frac{g-1}{2}$. If we observe that the condition $\deg
L>2g-\frac{g-1}{6}-2h^1(C,L)$ is equivalent to
$\cli(K_CL^{-1})<\frac{g-1}{6}$, then we meet an absurdity. Thus we
have $\deg N_2\le g-1$, and then Castelnuovo's genus bound produces
$\deg N_2\ge 3r(N_2)-2$. Note that the Castelnuovo number $\pi(d,r)$
has the property $\pi(d,r)\le\pi(d-2,r-1)$ for $d\ge 3r-2$ and $r\ge
3$, where $\pi(d,r)=\frac{m(m-1)}{2}(r-1)+m\epsilon,$
$d-1=m(r-1)+\epsilon,~~0\le \epsilon\le r-2$ (Lemma 6, \cite{KK1}).
Hence
$$
\pi(\deg N_2, r(N_2))\le \cdots\le \pi(\deg N_2-\deg G,
r(N_2)-\frac{\deg G}{2})\\
\le \pi (\deg N_1, r(N_1)),
$$
because of $2\le r(N_1)\le r(N_2)-\frac{\deg G}{2}$. Since
$r(N_1)\ge \frac{\deg N_1 }{4}$ and ${\deg N_1 }<\frac{g-1}{2}$, we
can induce a strict inequality $\pi(\deg N_1, r(N_1))<g$ as only the
number regardless of birational embedding from the proof of Lemma 9
in \cite{KK1}. It is absurd. Hence $m_2=3$, which yields $C_1\cong
C_2$.

Set $H_2$ be a hyperplane section of $C_2$. If $|H_2|$ on a smooth
model of $C_2$ is special, then $n:=r(N_2)\le \frac{\deg N_2}{6}$.
Thus the condition $\deg K_CL^{-1}<4(h^0(C,K_CL^{-1})-1)$ yields the
following inequalities:
$$
\frac{2\deg N_2}{3}\le\cli(N_2)\le \cli(N_1)\le \frac{\deg N_1}{2},
$$
which contradicts to $N_1\lneq N_2$. Accordingly $|H_2|$ is also
nonspecial.

Now we have $r(N_i)=\frac{\deg N_i}{3}-p, ~i=1,2$ where $p$ is the
genus of a smooth model of $C_1\cong C_2$. Therefore
$$
\frac{\deg N_1}{3}+2p=\cli(N_1)\ge \cli(N_2)=\frac{\deg N_2}{3}+2p
$$
which is a contradiction that $\deg N_1<\deg N_2$. This
contradiction comes from the assumption that $L$ is not normally
generated, thus the result follows.
\end{proof}

Using the above theorem, we obtain the following corollary under the
same assumption:

\begin{cor}
Let $C$ be a triple covering of genus $p$ curve $C'$ with
$C\stackrel{\phi}\rightarrow C'$ and $D$ a divisor on $C'$ with
$4p<\deg D< \frac{g-1}{6}-2p$. Then $K_C(-\phi^*D)$ becomes a very
ample line bundle which is normally generated.
\end{cor}

\begin{proof}
Set $d:=\deg D$ and $L:=K_C(-\phi^*D)$. Suppose $L$ is not base
point free, then there is a $P\in C$ such that
$|K_CL^{-1}(P)|=g^{r+1}_{3d+1}$. Note that $g^{r+1}_{3d+1}$ cannot
be composed with $\phi$ by degree reason. Therefore we have $g\le
6d+3p$ due to the Castelnuovo-Severi inequality. Hence it cannot
occur by the condition $d< \frac{g-1}{6}-2p$. Suppose $L$ is not
very ample, then there are $P,Q\in C$ such that
$|K_CL^{-1}(P+Q)|=g^{r+1}_{3d+2}$. By the same method as above, we
get a similar contradiction. Thus $L$ is very ample. The condition
$d<\frac{g-1}{6}-2p$ produces $\cli(K_CL^{-1})=d+2p<\frac{g-1}{6}$
since $\deg K_CL^{-1}=3d$ and $h^0(C,K_CL^{-1})=h^0(C', D)=d-p+1$.
Whence $\deg L>2g-\frac{g-1}{6}-2h^1(C,L)$ is satisfied. The
condition $4p<d$ induces $\deg K_CL^{-1}>4(h^0(C, K_CL^{-1})-1)$,
i.e., $\deg L>2g+2-4h^1(C,L)$. Consequently $L$ is normally
generated by Theorem \ref{3.5.6}.
\end{proof}

\begin{rmk}
In fact, we have similar result in \cite{KK1} for trigonal curve
$C$: $K_C(-rg^1_3)$ is normally generated if ~$3r<\frac{g}{2}-1$
$($\cite{KK}$)$. Thus our result could be considered as a
generalization which deals with triple covering under the some
condition.
\end{rmk}

\section{Application to projective surfaces}

Let $S\subseteq \mathbb P^r$ be a nondegenerate smooth surface and
$H$ a smooth hyperplane section of $S$. If $H$ is projectively
normal and $h^1(H, {\mathcal O}_H(2))=0$, then $q=h^1(S,{\mathcal
O}_S)=0, ~p_g=h^2(S,{\mathcal O}_S)=h^1(H,{\mathcal O}_H(1))$ and
$h^1(S,{\mathcal O}_S(t))=0$ for all nonnegative integer $t$
(\cite{Ak}, Lemma 2.1, Lemma 3.1). In this section, using our result
about the projective normality of smooth curves in section 2, we can
characterize smooth projective surfaces with the wider range of
degree and sectional genus. Recall the definition of $\Delta$-genus
given by $\Delta:=\deg S-r+1$.

\begin{thm}
\label{thm5.2} Let $S\subset \mathbb P^r$ be a nondegenerate smooth
surface of degree $2\Delta-e$ with $g(H)=\Delta+f$,
$\max\{\frac{e}{2}~, ~6e-\Delta\}<f-1<\frac{\Delta-2e-6}{3}$ for
some $e,f\in \mathbb Z_{\ge 1}$ and its general hyperplane section
$H$ is linearly normal. Then $S$ is projectively normal with
$~p_g=f$ and $-2f-e+2\le K_S^2\le \frac{(2f+e-2)^2}{2\Delta-e}$.
\end{thm}
\begin{proof}
From the linear normality of $H$, we get $h^0(H,{\mathcal O}_H(1))=
r$ and hence
\begin{eqnarray*}
h^1(H,{\mathcal O}_H(1))&=&-\deg{\mathcal O}_H(1)-1+g(H)+h^0(H,{\mathcal O}_H(1))\\
&=&-2\Delta+e-1+g(H)+h^0(H,{\mathcal O}_H(1))\\
&=& g(H)-\Delta=f
\end{eqnarray*}
Therefore we have $h^1(H,{\mathcal O}_H(1))>\deg \frac{K_H\otimes
\mathcal O_H(-1)}{4}$ since $f>\frac{e}{2}+1$ and $\deg {\mathcal
O}_H(1)=2\Delta-e=2g(H)-2-(2f+e-2)$. Thus ${\mathcal O}_H(1)$
satisfies $\deg {\mathcal O}_H(1)>2g(H)+2-4h^1(H,{\mathcal
O}_H(1))$. The condition $f-1> 6e-\Delta$ implies $\deg {\mathcal
O}_H(1)> 2g-\frac{g-1}{6}-2h^1(H,{\mathcal O}_H(1))$. Also the
condition $f-1<\frac{\Delta-2e-6}{3}$ yields $\deg {\mathcal
O}_H(1)>\frac{3g+3}{2}$. Hence ${\mathcal O}_H(1)$ is normally
generated by Theorem \ref{3.5.6}, and thus its general hyperplane
section $H$ is projectively normal since it is linearly normal.
Therefore $S$ is projectively normal with $q=0$,
$p_g=h^0(S,K_S)=h^1(H,{\mathcal O}_H(1))= f> 1$ since
$h^1(H,{\mathcal O}_H(2))=0$ from $\deg {\mathcal
O}_H(1)>\frac{3g+3}{2}$.

If we consider the adjunction formula $g(H)=\frac{K_S.H+H.H}{2}+1$
then $K_S.H=2f+e-2$. Since $|H+K_S|$ is ample and $p_g>0$, we get
$K_S.(H+K_S)\ge 0$. Therefore we obtain
$$K_S^2\ge H^2-(2g(H)-2)=(2\Delta-e)-(2(\Delta+f)-2)=-2f-e+2$$ by
Propositon 2.0.6 (iii) in \cite{Liv}. Thus $-2f-e+2\le K_S^2\le
\frac{(2f+e-2)^2}{2\Delta-e}$ by the Hodge index theorem
$K_S^2H^2\le (K_S. H)^2$. Hence the theorem is proved.
\end{proof}

Assume that $(2f+e-2)^2< 2\Delta-e$ in the above theorem, then we
have $-2f-e+2\le K_S^2\le 0$. Observe the cases for $K^2_S=-2f-e+2$
or $0$, then we obtain the following result by using similar method
in \cite{Ak}.

\begin{Prop}
\label{thm5.3} Let $S$ satisfy the conditions in Theorem
\ref{thm5.2}. Then $S$ is a minimal elliptic surface of Kodaira
dimension 1 if $K_S^2=0$ and $|K_S|$ has no fixed component. Also
$S$ is a surface blown up at $2f+e-2$ points on a $K3$ surface in
case $K_S^2=-2f-e+2$.

\end{Prop}
\begin{proof} Assume $|K_S|$ has no fixed component with
$K_S^2=0$. Then $S$ is minimal by adjunction formula and useful
remark III.5 in \cite{B}. Also the Kodaira dimension $\kappa$ of $S$
is at most one since $K_S^2\le 0$. Since $p_g>1$, $S$ is nonruled
and so $\kappa\ge 0$. If $\kappa=0$ then $p_g\le 1$ by Theorem
VIII.2 in \cite{B} and thus $\kappa$ must be 1. Hence by Proposition
IX.2 in \cite{B} there is a smooth curve $B$ and a surjective
morphism $p:S\to B$ whose generic fibre is an elliptic curve which
means that $S$ is a minimal elliptic surface of Kodaira dimension 1.

If $K_S^2=-2f-e+2$. Let $\phi_{H+K_S}=s\circ r$ be the
Remmert-Stein factorization of $\phi_{H+K_S}$ and $\hat{S}=r(S)$.
Then we can use Propositon 2.0.6 in \cite{Liv} as stated in the
proof of the previous theorem. And we obtain
$$H^2-K_S^2=(2\Delta-e)-(-2f-e+2)=2g(H)-2,$$ which yields $\hat{S}$ is
a minimal model and $K_{\hat{S}}= 0$, in other words, $\hat{S}$ is a
K3 surface by using Propositon 2.0.6 (iv-1) in \cite{Liv}. Also by
Propositon 2.0.6 (ii) in \cite{Liv}, $S$ is a surface blown up at
$2f+e-2$ points on a $K3$ surface $\hat{S}$ since
$\hat{d}-d=2g(H)-2-(2\Delta-e)=2f+e-2$.
\end{proof}

\pagestyle{myheadings} 

\end{document}